\documentclass[10pt]{amsart}

\usepackage{amsmath,amssymb,amsthm}
\usepackage[hidelinks]{hyperref}
\newtheorem{theorem}{Theorem}[section]
\newtheorem{proposition}[theorem]{Proposition}
\newtheorem{lemma}[theorem]{Lemma}
\newtheorem{corollary}[theorem]{Corollary}
\theoremstyle{definition}
\newtheorem{definition}[theorem]{Definition}
\theoremstyle{remark}
\newtheorem{remark}[theorem]{Remark}
\newcommand{\PP}{\mathbf P}
\newcommand{\CC}{\mathbf C}
\newcommand{\ZZ}{\mathbf Z}
\newcommand{\QQ}{\mathbf Q}
\newcommand{\OO}{\mathcal O}
\newcommand{\cC}{\mathcal C}

\DeclareMathOperator{\Br}{Br}
\DeclareMathOperator{\NS}{NS}
\DeclareMathOperator{\Pic}{Pic}
\DeclareMathOperator{\rk}{rk}

\DeclareMathOperator{\Gram}{Gram}

\title[Rationality of cubic fourfolds containing a plane]
{Rationality of cubic fourfolds containing a plane}
\author{Sumit K Singh}
\address{Department of Mathematics, Rutgers University, Piscataway, NJ 08854}
\email{ss4681@rutgers.edu}
\date{November 23, 2026}
\hypersetup{pdftitle={Rationality of cubic fourfolds containing a plane},
pdfauthor={Sumit K Singh}}

\begin{document}
\begin{abstract}
Let \(X\) be a smooth complex cubic fourfold containing a plane.
We prove that \(X\) is rational if and only if its full algebraic lattice
is admissible. The sufficiency argument is reduced to rank three.
For odd intersection parameter, Hassett's criterion gives a rational section of
the plane-projection bundle. In the remaining admissible cases we use
either an ineffective even theta characteristic of the discriminant 
curve with the same Clifford
class, or a rational section of the second projection of a birational Verra model.
The necessity direction uses
Gu\'er\'e's recent integral Hodge isometry of transcendental lattices and the
admissibility criterion of Addington and Thomas.
\end{abstract}
\maketitle
\tableofcontents

\section{Introduction}

The rationality problem for smooth cubic fourfolds has a long history.
Classical rationality constructions for cubic fourfolds were
studied by Fano \cite{F} and Tregub
\cite{T1,T2}, using quartic rational scrolls,
quintic del Pezzo surfaces, and pairs of disjoint planes.
For the rationality of Pfaffian cubic fourfolds, see
Beauville and Donagi \cite[Proposition~5(ii)]{B}.   

Hassett intiated the Hodge-theoretic study of smooth cubic fourfolds in \cite{Has00}.
He introduced the divisors $\mathcal C_d$ in the moduli
space $\mathcal C$, parametrizing smooth cubic fourfolds
whose algebraic lattice contains a primitive rank-two
sublattice of determinant $d$ containing the square of
the hyperplane class, and determined which discriminants
admit an associated polarized K3 surface.
We write \(A(X)=H^4(X,\mathbf Z)\cap H^{2,2}(X)\) for the full algebraic lattice 
of the cubic fourfold \(X\). Admissibility is recalled in §2.1.
The relation with Kuznetsov's categorical approach in \cite{Kuz}
was studied by Addington and Thomas \cite{AT},
who proved that the Hodge-theoretic and categorical conditions agree
generically on the admissible Hassett divisors.
Russo and Staglian\`o \cite{RS1,RS2} proved that every smooth cubic 
fourfold in the Hassett divisors \(\mathcal C_{26}\), \(\mathcal C_{38}\), 
and \(\mathcal C_{42}\) is rational.
Numerous examples have also been found in 
codimension two; see \cite{Has99,AH,YW,Has24}.

Iritani \cite{It} established a decomposition theorem for
quantum cohomology under blow-up.
Building on this result, Katzarkov, Kontsevich, Pantev, and Yu
\cite{KK} introduced Hodge atoms.
Their behavior under blow-ups gives new birational obstructions;
in particular, these authors proved that a very general cubic
fourfold is irrational. These authors also announce the irrationality 
of a very general cubic fourfold containing a plane; see \cite[Example 6.17]{KK}.
Gu\'er\'e \cite{GuereCubic} used this framework to obtain a
necessary K3 condition for the rationality of cubic fourfolds at
the level of rational Hodge structures.

The present paper concerns rationality for cubic fourfolds 
that contain a plane. The projection from the plane relates 
their rationality to that of quadric surface bundles.
The description of quadric surface bundles by discriminant covers and
Clifford algebras provides a way to study their birational geometry.
Auel, Bernardara, and Bolognesi \cite{ABB} studied the behavior of
the even Clifford algebra under quadric reduction.
Auel, Bernardara, Bolognesi, and V\'arilly-Alvarado \cite{ABBV}
exhibited rational cubic fourfolds containing a plane for which
the Clifford class of the plane projection is nonzero.
These examples show that rationality need not arise from a section
of the original projection.
Bini, Kapustka, and Kapustka \cite{BKK} studied the
birational models of quadric surface bundles and proved a
rationality result for bundles associated with even theta
characteristics on smooth plane sextics
\cite[Corollary~4.11]{BKK}.

The main result is as follows:

\begin{theorem}\label{thm:main}
Let \(X\subset\PP^5\) be a smooth complex cubic fourfold containing
a plane. Then
\[X\text{ is rational}\quad\Longleftrightarrow\quad
 A(X)\text{ is admissible}.\]
\end{theorem}

\noindent\textbf{Outline of Proof}
The sufficiency(\(\Leftarrow\)) argument the rank-three lattices
\begin{equation}\label{eq:M}
 M_{\tau,n}=
 \begin{pmatrix}
 3&1&0\\
 1&3&\tau\\
 0&\tau&2n
 \end{pmatrix},
 \qquad
 \det M_{\tau,n}=16n-3\tau^2.
\end{equation}
The ordered basis is \((H,P,R)\), and the parameters satisfy
\[
 0\leq\tau\leq4,\qquad n\geq2,\qquad
 (\tau,n)\notin\{(3,2),(4,2),(4,3)\}.
\]
Yang and Yu classified the rank-three loci among cubic fourfolds containing a plane
and determined their admissibility \cite[Corollary 8.14]{YY}.

For \(\tau=1,3\), the class \(Q=H-P\) of the quadric surface residual to \(P\) 
pairs oddly with the third generator \(R,\) and Hassett's theorem applies 
\cite[Theorem 3.1]{Has99}.

For \(\tau=0,2,4\) we need additional arguments. 
For \(\tau=0,4\) and odd \(n\), the calculation of the vanishing Brauer class in
\cite{GvGbrauer} supplies an even theta characteristic encoding the same 
Clifford class as the plane-projection bundle of \(X\).
A Picard-lattice calculation shows that this theta characteristic is ineffective.
Its symmetric minimal resolution defines a net of quadrics in \(\PP^5\). 
Choose a point in their common base locus. Quadric reduction 
at this point gives a quadric surface bundle \(Z\to\mathbf P^2\) whose total space is 
rational. Its discriminant cover and Clifford class agree with those of the 
original plane-projection bundle of \(X\). 
The generic quadric surfaces are therefore 
isomorphic over \(\mathbf C(\mathbf P^2)\), giving an isomorphism 
\(\mathbf C(Z)\cong\mathbf C(X)\); see \cite{BKK}.
 
For a very general member \(X\) with \(\tau=2\) and odd \(n\), the construction 
in \cite{GvGVerra} gives a smooth Verra fourfold birational to \(X\). 
After labeling its projections, 
the first quadric surface bundle is birational over \(\mathbf P^2\) to the 
original plane-projection bundle of \(X\). The resulting integral Hodge isometry of 
transcendental lattices, together with a calculation of discriminant bilinear 
pairings, forces a rational section of the second projection. This proves 
the rationality 
of the very general member. The specialization theorem of Kontsevich and
Tschinkel \cite{KT} extends these conclusions to every smooth member
of the corresponding marked loci. The saturation argument then treats
an admissible lattice of arbitrary rank.

The necessity(\(\Rightarrow\)) argument combines the integral decomposition of transcendental 
lattices in \cite[Propositions 25–26]{GuereIntegral} with the analysis of 
surface centers in 
\cite[proof of Theorem 63]{GuereCubic}. Proposition 3.13 explains how these results 
yield an integral Hodge isometry \(T(X)\simeq T(S)(-1)\) for a projective 
K3 surface \(S\). The criterion in \cite{AT} then gives admissibility.
Thus, the necessity direction, and
hence the equivalence in Theorem~\ref{thm:main}, uses the recent
preprint \cite{GuereIntegral}.

This discussion culminates in the following result:

\begin{corollary}\label{cor:rank-three}
Suppose that \(A(X)=M_{\tau,n}\), with the parameters in the range
above. Then \(X\) is rational if and only if
\[\tau\in\{1,3\}\quad\text{or}\quad
\bigl(\tau\in\{0,2,4\}\text{ and }n\text{ is odd}\bigr).\]
For each admissible pair, every smooth cubic fourfold with a primitive
\(M_{\tau,n}\)-marking is rational.
\end{corollary}

\medskip
\noindent\textbf{Outline.}
Section~2 recalls the background on algebraic and transcendental
lattices, Clifford classes of quadric surface bundles, theta
characteristics, and Verra fourfolds. It also explains the results
on transcendental lattices and birational maps used in the
necessity argument. In Section~3.1, we recall the classification
and admissibility of the rank-three lattices in $\mathcal C_8$.
Sections~3.2 and~3.3 establish the rationality results using theta
characteristics and Verra fourfolds, respectively. Section~3.4
extends these results by smooth specialization and treats
admissible algebraic lattices of higher rank. Finally,
Section~3.5 proves the necessity of admissibility and completes
the proofs of Theorem~1.1 and Corollary~1.2.

\medskip
\noindent\textbf{Acknowledgement}
I thank Lev Borisov and Lisa Marquand for giving valuable ideas and insightful discussions. 

\medskip
\noindent\textbf{AI Disclosure}
The author acknowledges the use of AI in the research process and the preparation of this manuscript. The author assumes full responsibility for all content.

\section{Preliminaries}\label{sec:preliminaries}

\subsection{Cubic fourfolds and lattices}
Let \(X\subset\mathbf P^5\) be a smooth cubic fourfold over \(\mathbf C\),
and write
\[H=h_X^2,\qquad A(X)=H^4(X,\mathbf Z)\cap H^{2,2}(X),\qquad T(X)=A(X)^\perp\subset H^4(X,\mathbf Z),\]
where \(h_X=c_1(\mathcal O_X(1))\). The Hassett divisor
\(\mathcal C_8\) parametrizes cubic fourfolds containing a plane.
If \(P\subset X\) is such a plane, the projection from \(P\) gives a
quadric surface bundle
\[q_X:\operatorname{Bl}_P X\longrightarrow\mathbf P^2.\]
A rational section makes the generic quadric rational over
\(\mathbf C(\mathbf P^2)\), and hence makes \(X\) rational.

The intersection form on \(H^4(X,\ZZ)\) is unimodular of signature
\((21,2)\), and \(H^2=3\). The lattice \(A(X)\) is positive definite
and primitive. An equality \(A(X)=M\) will mean an identification
with the full algebraic lattice, including the distinguished class
\(H\). By contrast, an \(M\)-marking is a primitive isometric embedding
\(M\hookrightarrow A(X)\) taking its distinguished class to \(H\).
The latter permits \(\rk A(X)>\rk M\).

A discriminant-\(d\) marking is a primitive rank-two sublattice
\(K_d\subset A(X)\) containing \(H\), with determinant \(d\).
The corresponding Hassett divisor \(\cC_d\) is nonempty precisely
when \(d>6\) and \(d\equiv0,2\pmod6\).
Such a discriminant is called admissible if
\begin{equation}\label{eq:admissible-d}
  4\nmid d,\qquad 9\nmid d,\qquad
p\nmid d\quad\text{for every odd prime }p\equiv2\pmod3.
\end{equation}

We call \(A(X)\) admissible if it contains a primitive rank-two sublattice containing \(H\) whose discriminant is admissible.
Condition
\eqref{eq:admissible-d} is the condition for an associated polarized
K3 surface; see \cite[Theorems~1.0.1--1.0.2 and~5.1.3]{Has00}.

For a projective K3 surface \(S\), put
\[
 T(S)=\NS(S)^\perp\subset H^2(S,\ZZ).
\]
We will also use the following equivalent form of admissibility:
\begin{equation}\label{eq:associated-K3}
 A(X)\text{ is admissible}
 \quad\Longleftrightarrow\quad
 T(X)\simeq T(S)(-1)
\end{equation}
for some projective K3 surface \(S\), where the isometry is integral
and respects the Hodge structures.

 The forward implication follows
by restricting the Hodge isometry defining an associated \(\operatorname{K3}\) surface to the transcendental lattices. For the converse, the
extension argument in \cite[Section~8, proof of Theorem~1.3]{AT}
extends the isometry of  transcendental lattices to the topological Mukai lattice
of \(S\) and \(\operatorname{Ku}(X).\) Kuznetsov component of a cubic fourfold  
\(\operatorname{Ku}(X)\) is defined in \cite{Kuz}.
The algebraic Mukai lattice then contains a hyperbolic plane, and
\cite[Theorem~3.1]{AT} provides an admissible marking. There is no
rank-three assumption in this equivalence.

All lattices below are integral and nondegenerate. A sublattice is
primitive if the quotient is torsion free. For a lattice \(M\), write
\[
 M^\vee=\{x\in M\otimes\QQ:(x,M)\subset\ZZ\},\qquad
 D(M)=M^\vee/M,
\]
and denote its discriminant bilinear pairing by
\[
 b_M:D(M)\times D(M)\longrightarrow\QQ/\ZZ,
 \qquad b_M([x],[y])=(x,y)\bmod\ZZ.
\]
Thus \(|D(M)|=|\det M|\). Define the subgroup 
\(D(M)[2]=\{\xi\in D(M):2\xi=0\}.\)
For a sublattice \(M\subset L\), define the saturation or primitive closure as
\[M_L^{\mathrm{sat}}=(M\otimes_{\mathbf Z}\mathbf Q)\cap L.\]

We use signed determinants and write
\[
 U=\begin{pmatrix}0&1\\1&0\end{pmatrix}.
\]
The notation \(M(a)\) multiplies the pairing by \(a\). 
For Hodge lattices, \((-1)\) denotes the Tate twist, which shifts bidegrees by 
\((1,1)\); we equip the twisted lattice with the negative of the original pairing.

\begin{lemma}\label{lem:discriminant}
Let \(M\) be a primitive nondegenerate sublattice of a unimodular
lattice \(L\), and put \(T=M^\perp\). Then
\[(D(M),b_M)\simeq(D(T),-b_T).\]
\end{lemma}
\begin{proof}
The quotient \(\Gamma=L/(M\oplus T)\) embeds into
\(D(M)\oplus D(T)\). Primitivity makes both projections injective.
Unimodularity gives
\[
 |\Gamma|^2=|D(M)|\,|D(T)|.
\]
Since both projections are injective, \(|\Gamma|\leq |D(M)|\) and 
\(|\Gamma|\leq |D(T)|\); together with the formula for the displayed product, 
these inequalities force equality.
Hence \(|\Gamma|=|D(M)|=|D(T)|\), and \(\Gamma\) is the graph of
an isomorphism \(D(M)\to D(T)\). The integrality of the pairing on \(L\)
shows that the two discriminant pairings sum to zero on this graph.
\end{proof}

This statement applies to odd ambient lattices. In particular, an
integral isometry of transcendental lattices in unimodular middle
cohomology identifies the discriminant bilinear pairings of their
primitive algebraic orthogonal complements.

\subsection{Clifford classes of quadric surface bundles}

Suppose that $X$ contains a plane $P$. We also write $P$ for
its class in $A(X)$. The sublattice generated by $H$ and $P$
is primitive and has Gram matrix
\[
\operatorname{Gram}(H,P)=
\begin{pmatrix}
3&1\\
1&3
\end{pmatrix}.
\]
Projection from $P$ induces a quadric surface bundle
\[
q_X\colon \widetilde X=\operatorname{Bl}_P X
\longrightarrow \mathbf P^2.
\]
The class in $A(X)$ of a quadric surface residual to $P$ is
\[
Q=H-P,
\qquad
Q\cdot H=2,
\qquad
Q\cdot P=-2.
\]
Thus $Q$ is the pushforward to $X$ of the class of a fiber
of $q_X$.

Assume that the discriminant sextic $C$ of $q_X$ is smooth.
This holds, in particular, when no second plane in $X$
meets $P$; see \cite[Section~1]{ABBV}. The double cover
\[\pi\colon B\longrightarrow \mathbf P^2\]
branched along $C$ is a K3 surface. We put
\[h=c_1\bigl(\pi^*\mathcal O_{\mathbf P^2}(1)\bigr),\qquad h^2=2.\]

Here simple degeneration means that each singular geometric fiber has co-rank one.
More generally, let $q\colon W\to\mathbf P^2$ be a quadric
surface bundle with simple degeneration along a smooth
discriminant sextic $C$, and let $B$ be its discriminant
double cover. The even Clifford algebra determines an
Azumaya algebra on $B$. We denote its Brauer class by
\[
\alpha_W\in\operatorname{Br}(B)[2]
\]
and call it the Clifford class of $q$. The relative Fano
scheme of lines has Stein factorization
\[
F_1(W/\mathbf P^2)\longrightarrow B
\longrightarrow\mathbf P^2.
\]
The first morphism is a Severi--Brauer conic with Brauer
class $\alpha_W$; see \cite[Appendix~B,
Propositions~B.2--B.3]{ABB}.

Put $K=\mathbf C(\mathbf P^2)$ and $L=\mathbf C(B)$.
If $\varphi_\eta$ is a quadratic form defining the generic
fiber of $q$, then its even Clifford algebra is a quaternion
algebra over its center $L$, and
\[
\alpha_W|_L=[C_0(\varphi_\eta)]
\quad\text{in }\operatorname{Br}(L).
\]
This follows from compatibility with base change
\cite[Proposition~1.18 and Proposition~A.2(5)]{ABB}.
Moreover,
\[
\alpha_W=0
\quad\Longleftrightarrow\quad
q\text{ has a rational section};
\]
see \cite[Theorem~1.33]{ABB}. For the plane projection, we write
$\alpha_X=\alpha_{\widetilde X}$.

We also use the following form of Hassett's criterion.
Its statement does not require a smooth discriminant.

\begin{proposition}\label{death}
Let $q\colon W\to\mathbf P^2$ be a quadric surface bundle
with smooth and projective \(W\) and with generic smooth
fiber. Let $f\in H^4(W,\mathbf Z)$ be the class of a smooth
fiber. If
\[
z\in H^4(W,\mathbf Z)\cap H^{2,2}(W)
\qquad\text{and}\qquad
z\cdot f\text{ is odd},
\]
then $q$ has a rational section and $W$ is rational.
\end{proposition}

\begin{proof}
We follow \cite[Propositions~2.1 and~2.3]{Has99}.
Retain the components of the relative Fano scheme of lines
$F$ that dominate $\mathbf P^2$, and choose a smooth projective
resolution $\widetilde F$ of their union. The incidence
correspondence sends $z$ to an integral $(1,1)$-class on
$\widetilde F$. By the Lefschetz $(1,1)$-theorem, this is
the first Chern class of a line bundle $\mathcal L$.

Let $F_\eta$ be the Fano scheme of the generic quadric.
Its two geometric components correspond to the two
rulings. The restriction $\mathcal L_\eta$ has degree
$m=z\cdot f$ in each component. Hence
\[
\mathcal M=
\mathcal L_\eta\otimes
\omega_{F_\eta/K}^{\otimes (m-1)/2}
\]
has degree one on each component.

 After extending scalars to $\overline K$, the two components
are copies of $\mathbf P^1$, and $\mathcal M$ restricts to
$\mathcal O_{\mathbf P^1}(1)$ on each component. By cohomology
and base change, and since $K$ is infinite, we may choose a
section of $\mathcal M$ defined over $K$ whose restriction
to each geometric component is nonzero.
Its zero divisor gives a Galois-invariant pair of lines in opposite
rulings. Their intersection is a $K$-point of the generic
quadric and therefore gives a rational section of $q$.
The Projection from this point identifies its function field
with $K(u,v)$. Since $K=\mathbf C(\mathbf P^2)$, the
fourfold $W$ is rational.
\end{proof}

For the plane projection, the criterion takes the form
\begin{equation}\label{life}
q_X\text{ has a rational section}
\quad\Longleftrightarrow\quad
D\cdot Q\text{ is odd for some }D\in A(X);
\end{equation}
see \cite[Theorem~3.1]{Has99} and \cite[proof of Proposition~2]{ABBV}.

We now return to a smooth discriminant sextic $C$ and its
double cover $B$. Adjunction gives
\[
\omega_C\simeq\mathcal O_C(3),
\qquad g(C)=10.
\]

\begin{definition}
 A theta characteristic on $C$ is a line bundle $\theta$
such that $\theta^{\otimes 2}\simeq\omega_C$.
\end{definition}

Its parity
is the parity of $h^0(C,\theta)$, and it is ineffective
when $h^0(C,\theta)=0$.
In a smooth proper family of complex curves, a theta characteristic on one fiber extends locally in the analytic topology, or étale locally in the algebraic setting. For a family of theta characteristics \(\theta_t\), the parity
\[
h^0(C_t,\theta_t)\pmod 2
\]
is locally constant. In particular, an even theta characteristic remains even under smooth specialization; see \cite[\S ~1]{M}.

Let $[K_C]$ denote the canonical class in
$\operatorname{Pic}(C)$. There is an exact sequence
\begin{equation}\label{Sumit}
0\longrightarrow
\frac{\operatorname{Pic}(B)}
     {\mathbf Z h+2\operatorname{Pic}(B)}
\longrightarrow
\left(
\frac{\operatorname{Pic}(C)}{\mathbf Z[K_C]}
\right)[2]
\xrightarrow{\ \mathfrak b_B\ }
\operatorname{Br}(B)[2]
\longrightarrow 0;
\end{equation}
see \cite[Theorem~1.1]{IOOV}. The homomorphism $\mathfrak b_B$
agrees with the construction using Clifford algebras
\cite[Theorem~6.1]{IOOV}.

A theta characteristic defines an element of the middle
group because $2[\theta]=[K_C]$. We write
$\mathfrak b_B(\theta)$ for the image of this element.
Thus, this notation refers to the class of $\theta$ in
the quotient, rather than to a two-torsion line bundle
in $\operatorname{Pic}(C)$. The plane projection determines
an odd theta characteristic $\theta_X$ satisfying
\[
h^0(C,\theta_X)=1,
\qquad
\mathfrak b_B(\theta_X)=\alpha_X;
\]
see \cite[Section~3.1]{GvGbrauer} and \cite[Theorem~6.1]{IOOV}.

For a complex K3 surface $B$, there is a natural
identification
\[\operatorname{Br}(B)[2]\simeq
\operatorname{Hom}\bigl(T(B),\mathbf Z/2\mathbf Z\bigr).\]

We identify $\mathbb Z/2\mathbb Z$ with
$(\frac12\mathbb Z)/\mathbb Z$ by sending $1$ to
$\frac12+\mathbb Z$.
For $\alpha\in\operatorname{Br}(B)[2]$, a
$B$-field representative of $\alpha$ is a class
\[
B_\alpha\in\frac12 H^2(B,\mathbb Z)
\]
such that
\[
\alpha(x)=B_\alpha\cdot x\pmod{\mathbb Z}
\qquad\text{for every }x\in T(B).
\]
Every two-torsion Brauer class admits such a representative;
see \cite[Sections~1.1--1.2]{GvGbrauer}.

We will also use specialization of Brauer classes in smooth families of K3 
surfaces. After choosing a local marking, suppose that \(B_t\) specializes to 
\(B_0\) with
\[
\operatorname{NS}(B_t)\subseteq\operatorname{NS}(B_0),
\qquad T(B_0)\subseteq T(B_t).
\]
Under the preceding description of the Brauer group, specialization is the restriction map
\[
\operatorname{sp}:\operatorname{Br}(B_t)[2]
\longrightarrow\operatorname{Br}(B_0)[2],
\qquad
\alpha\longmapsto\alpha|_{T(B_0)}.
\]
In particular, if \(\delta\) has a \(B\)-field representative \(\kappa/2\), where \(\kappa\) becomes algebraic on \(B_0\), then
\[
\operatorname{sp}(\delta)(x)
=\frac{\kappa\cdot x}{2}=0
\quad\text{in }\mathbf Q/\mathbf Z
\qquad (x\in T(B_0)).
\]
Consequently, \(\operatorname{sp}(\alpha+\delta)=\operatorname{sp}(\alpha)\);
see \cite[Definition~2.1.1 and Proposition~2.1.2]{GvGbrauer}.

For $\alpha\in\operatorname{Br}(B)[2]$, put
\[T(B,\alpha)=\ker\bigl(\alpha\colon T(B)\to\mathbf Z/2\mathbf Z\bigr).\]

The plane projection induces an integral Hodge isometry
\begin{equation}
T(X)\simeq T(B,\alpha_X)(-1),
\end{equation}
and
\[23-\operatorname{rk}A(X)
=\operatorname{rk}T(X)
=\operatorname{rk}T(B)
=22-\rho(B).\]
see \cite[Proposition~3.2.1]{GvGbrauer}. No assumption on the Picard
rank of $B$ is required here.

\begin{proposition}\label{quad}
Let $W_1,W_2\to\mathbf P^2$ be quadric surface bundles
with smooth generic fibers. Suppose that the discriminant
quadratic étale $K$-algebras of their generic fibers are
identified, and that their even Clifford algebras have
the same Brauer class under this identification.
Then the generic fibers are isomorphic over $K$, and
$W_1$ and $W_2$ are birational over $\mathbf P^2$.
\end{proposition}

\begin{proof}
The final assertion of \cite[Proposition~2.4]{BKK} identifies
the generic quadric surfaces over $K$. Their function
fields are therefore isomorphic as extensions of $K$,
which gives the required birational map over
$\mathbf P^2$; see also \cite[Corollary~2.5]{BKK}.
\end{proof}

Finally, consider a quadric bundle of even relative
dimension with simple degeneration along a smooth
divisor. Quadric reduction along a regular isotropic
subbundle preserves the discriminant cover and the
Brauer class of the even Clifford algebra
\cite[Proposition~3.2]{BKK}. Here, regular means that the
isotropic subbundle meets the radical trivially in
every geometric fiber. For an isotropic line subbundle,
the corresponding section therefore meets the smooth
locus of every fiber. This is the condition called
smoothness in \cite[Section~1.3]{ABB}; the preservation of the
Brauer class follows from \cite[Theorem~1.27]{ABB}.

\subsection{Verra fourfolds}
A smooth Verra fourfold is a double cover
\[ V\longrightarrow\PP^2\times\PP^2\]
branched along a smooth divisor of bidegree \((2,2)\). Its two
projections \(p_1,p_2\) are quadric surface bundles. Put
\[h_i=c_1(p_i^*\OO_{\PP^2}(1)),\qquad f_i=h_i^2,\qquad e=h_1h_2.\]
Then
\begin{equation}\label{eq:Verra-lattice}
 \langle f_1,f_2,e\rangle\simeq U(2)\oplus\langle2\rangle
\end{equation}
is primitive in \(H^4(V,\ZZ)\). This cohomology is torsion free and
unimodular, of rank \(24\) and signature \((21,3)\)
\cite[Section~1.2 and Proposition~1.3(a)]{GvGVerra}.
Write
\[A(V)=H^4(V,\ZZ)\cap H^{2,2}(V),\qquad T(V)=A(V)^\perp.\]
Suppose that the plane-projection bundle of \(X\) and one projection of \(V\) have an identified smooth discriminant double cover \(B\), and that their Clifford classes agree under this identification. Write \(\alpha\in\operatorname{Br}(B)[2]\) for the common class. The integral identifications for these two quadric surface bundles give
the following:
\[T(X)\simeq T(B,\alpha)(-1)\simeq T(V).\]
Thus, \(T(X)\) and \(T(V)\) are integrally Hodge isometric; see 
\cite[Proposition 3.2.1]{GvGbrauer} and \cite[Proposition 1.3(b) and §1.7]{GvGVerra}.

\subsection{Transcendental lattices and birational maps}
\label{subsec:transcendental-birational}

For a smooth projective variety \(W\), write
\[
H^j(W,\mathbf Z)_{\mathrm{free}}
=
H^j(W,\mathbf Z)/\operatorname{torsion}.
\]
If \(Y\) is a smooth projective surface, let
\(\operatorname{NS}(Y)_{\mathrm{free}}\) denote the image of its
N\'eron--Severi group in \(H^2(Y,\mathbf Z)_{\mathrm{free}}\), and put
\[
T(Y)
=
\operatorname{NS}(Y)_{\mathrm{free}}^\perp
\subset H^2(Y,\mathbf Z)_{\mathrm{free}}.
\]
The orthogonal complement is taken with respect to the intersection
form. We equip \(T(Y)\) with its induced integral Hodge structure.
In the notation \(T(Y)(-1)\), the Tate twist raises both Hodge indices
by one; the lattice pairing is multiplied by \(-1\).

Let \(\pi:\widetilde Y\to Y\) be the blow-up of a point, with exceptional
curve \(E\). The integral blow-up formula gives
\[
H^2(\widetilde Y,\mathbf Z)_{\mathrm{free}}
=
\pi^*H^2(Y,\mathbf Z)_{\mathrm{free}}
\mathbin{\perp}\mathbf Z[E],
\qquad E^2=-1.
\]
Moreover,
\[
\operatorname{NS}(\widetilde Y)_{\mathrm{free}}
=
\pi^*\operatorname{NS}(Y)_{\mathrm{free}}
\oplus\mathbf Z[E].
\]
Taking orthogonal complements therefore gives an integral Hodge
isometry
\[
\pi^*:T(Y)\xrightarrow{\sim}T(\widetilde Y).
\]
Consequently, birational smooth projective surfaces have isometric
transcendental lattices.

For fourfolds, surface centers must also be considered.
By the weak factorization theorem
\cite[Theorem~0.1.1]{AKM}, a birational map between smooth projective
fourfolds admits a factorization
\[
W=W_0\longleftrightarrow W_1
\longleftrightarrow\cdots
\longleftrightarrow W_m=W',
\]
where each step is a blow-up along a smooth center or its inverse,
and every \(W_i\) is smooth and projective. After omitting trivial
steps, the centers have dimension at most two. If
\(\widetilde W=\operatorname{Bl}_Y(W)\) with \(Y\) a smooth surface,
then
\[
H^4(\widetilde W,\mathbf Z)_{\mathrm{free}}
\simeq
H^4(W,\mathbf Z)_{\mathrm{free}}
\mathbin{\perp}
H^2(Y,\mathbf Z)_{\mathrm{free}}(-1).
\]
Here the second summand is induced by \(j_*p^*\), where
\(p:E\to Y\) is the exceptional divisor and
\(j:E\hookrightarrow\widetilde W\). Its pairing is the negative
intersection form of \(Y\). Centers of dimension zero or one add
only algebraic classes in degree four. Thus weak factorization
alone does not identify the transcendental lattices of birational
fourfolds: surface contributions can occur on both sides.

We recall the additional input used in Section~3.5.
Let \(\star_{\mathbf t}\) denote the big quantum product and let
\(\kappa_{\mathbf t}=\operatorname{Eu}_{\mathbf t}
\star_{\mathbf t}(-)\) be quantum multiplication by the Euler field;
see \cite[Section~1.1]{GuereIntegral}. Evaluation maps are understood
in the completed, graded coefficient rings of
\cite[Definition~19 and Proposition~20]{GuereCubic}.
After extending coefficients, set
\[
E^W_{\mathrm{ev},\alpha}
=
\ker\!\left(
\mathrm{ev}(\kappa_{\mathbf t})-\alpha b^2\operatorname{id}
\right)^N,
\qquad N\gg0,
\]
where \(b\) is the degree-one grading parameter.
For evaluations with maximal spectrum cardinality, these
generalized eigenspaces are the coarse Hodge atoms
\cite[Introduction and Definition~26]{GuereCubic}.

The distinction between generic and arbitrary evaluations is
essential. In \cite[Definition~33]{GuereCubic}, the generic
condition is imposed on a dense set of evaluations with maximal
spectrum cardinality, whereas its strong version is imposed on all
evaluation maps. \cite[Corollary~47]{GuereCubic} requires the strong condition on
every center. 
For Property \(\heartsuit\), defined in \cite[Definitions 30 and 33]{GuereCubic}, Proposition 60 in \cite{GuereCubic} identifies the generic and strong conditions 
in dimension at most two, over number fields containing \(\mathbb Q(i)\).
It does not assert that every surface
satisfies the condition.

The integral condition is separate. For a rational cubic fourfold
\(X\), \cite[Proposition~25]{GuereIntegral} gives surface centers
\(Y_1,\ldots,Y_r\) in a weak factorization to \(\mathbf P^4\) and
an orthogonal integral Hodge isometry
\[
T(X)\simeq\bigoplus_{i=1}^r T(Y_i)(-1).
\]
The even transcendental cohomology in that proposition equals
\(T(X)\), since the remaining even cohomology of a cubic fourfold
is algebraic. In \cite[Proposition~26 and Section~3.3]{GuereIntegral} of the same paper
identify the integral transcendental summand of an atom with that
of its corresponding surface center. This compatibility is needed
to use the same center in the lattice decomposition and in the
surface analysis of
\cite[Propositions~59--60 and proof of Theorem~63]{GuereCubic}.
An isomorphism after tensoring with \(\mathbf Q\) would not suffice.

Finally, if the resulting smooth minimal surface \(S\) satisfies
\[
c_1(K_S)=0\ \text{in }H^2(S,\mathbf Q),
\qquad
h^1(S,\mathcal O_S)=0,
\qquad
h^2(S,\mathcal O_S)=1,
\]
then \(S\) is a K3 surface. Indeed, a nonzero canonical section has
an effective zero divisor numerically equivalent to zero. Such a
divisor must vanish, since every nonzero effective divisor meets
an ample class positively. Hence \(K_S\simeq\mathcal O_S\), and
the vanishing of \(h^1(S,\mathcal O_S)\) gives the assertion.
The preceding point blow-up calculation then identifies the
transcendental lattice of the surface center with \(T(S)\).

\section{Rationality in \(\cC_8\)}\label{sec:rationality}

\subsection{Admissible rank-three lattices in \texorpdfstring{\(\cC_8\)}{C8}}
Let \(\cC_{\tau,n}\) denote the locus of smooth cubics with a
primitive \(M_{\tau,n}\)-marking. The allowed pairs are those in
Corollary~\ref{cor:rank-three}. Each such locus is nonempty,
irreducible, and of dimension \(18\); its very general member has
full algebraic lattice \(M_{\tau,n}\)
\cite[Corollary~8.14]{YY}.

\begin{proposition}[Yang--Yu]\label{prop:parity}
The lattice \(M_{\tau,n}\) is admissible if and only if
\[
 \tau=1,3\quad\text{or}\quad \tau=0,2,4\text{ with }n\text{ odd}.
\]
\end{proposition}
\begin{proof}
This is \cite[Corollary~8.14]{YY}, in the form of
\cite[Proposition~5.1.3]{GvGbrauer}.
\end{proof}

Let \(H,P,R\) be the standard basis of \(M_{\tau,n}\).
The class of the quadric surface residual to \(P\) satisfies
\begin{equation}\label{eq:Q-pairings}
 Q\cdot H=2,\qquad Q\cdot P=-2,\qquad Q\cdot R=-\tau.
\end{equation}
Thus \(\tau=1,3\) gives a rational section by
\eqref{life}. For even \(\tau\), the plane-projection
bundle of a cubic with full algebraic lattice \(M_{\tau,n}\) has no rational
section. We treat its total space by another birational model.

\subsection{Rationality through theta characteristics}
The following lemma also follows from the classification of symmetric resolutions 
in \cite[Theorem 5.1, Corollary 5.2 and Sections 7.4–7.5]{IOOV}. 
We include a direct proof.

\begin{lemma}\label{lem:ineffective}
Let \(C\subset\PP^2\) be a smooth plane sextic, and let
\(\pi:B\to\PP^2\) be the double cover branched along \(C\).
Put \(h=c_1(\pi^*\OO_{\PP^2}(1))\), and suppose that
\(h\cdot\NS(B)\subset2\ZZ\). Then every theta characteristic
\(\theta\) on \(C\) has \(h^0(C,\theta)\leq1\).
In particular, every even theta characteristic is ineffective.
\end{lemma}
\begin{proof}
Suppose that \(s,t\in H^0(C,\theta)\) are independent. Their products
\(s^2,st,t^2\) are sections of \(K_C=\OO_C(3)\).
Restriction identifies plane cubic polynomials with
\(H^0(C,\OO_C(3))\). Denote the corresponding cubics by \(F,G,J\):
\[
 F|_C=s^2,\qquad G|_C=st,\qquad J|_C=t^2.
\]
If \(f_C\) defines \(C\), then
\[
 FJ-G^2=c f_C
\]
for a constant \(c\in\CC\).

Suppose first that \(c\ne0\). Consider
\[
 W=\{Fu^2+2Guv+Jv^2=0\}\subset\PP^1_{[u:v]}\times\PP^2.
\]
The cubics \(F,G,J\) have no common zero, since such a zero would
make \(f_C\) and all its first derivatives vanish. Thus projection
to \(\PP^2\) is a finite flat double cover. Its branch curve is
\(C\), so \(W\simeq B\). A fiber \(E\) of the other projection is
a plane cubic and has \(h\cdot E=3\), a contradiction.

Suppose now that \(c=0\). Unique factorization gives
\[
 F=AU^2,\qquad G=AUV,\qquad J=AV^2,
\]
with \(U,V\) coprime. They cannot be constant, since \(s,t\) are
independent. Degree considerations imply that \(A,U,V\) are linear.
Taking zero divisors on \(C\) gives
\[
 2\operatorname{div}(s)=\operatorname{div}(A|_C)
                         +2\operatorname{div}(U|_C).
\]
Hence the line \(\ell=\{A=0\}\) meets \(C\) with even multiplicity
at every point. The restriction \(f_C|_\ell\) is a square. Its inverse
image in \(B\) therefore has two components, each of \(h\)-degree
one, again a contradiction. This proves \(h^0(C,\theta)\leq1\).
An even theta characteristic must consequently have no sections.
\end{proof}

The following is a consequence of Beauville’s determinantal description 
\cite[Proposition 4.2(a)]{Beauville}. We recall the construction and verify 
smoothness of the base locus.

\begin{lemma}\label{lem:net}
Let \(C\subset\mathbf P^2\) be a smooth plane sextic, with inclusion \(i\), and let \(\beta\) be an ineffective theta characteristic on \(C\).
Its symmetric minimal resolution 
defines a net of quadrics in \(\PP^5\) with smooth base locus.
\end{lemma}

\begin{proof}
By \cite[Proposition~4.2(a)]{Beauville}, there are symmetric
\(6\times6\) matrices \(M_0,M_1,M_2\) such that
\[
 M(\lambda)=\lambda_0M_0+\lambda_1M_1+\lambda_2M_2
\]
fits into the exact sequence
\[
 0\longrightarrow\OO_{\PP^2}(-2)^6
 \xrightarrow{\ M(\lambda)\ }\OO_{\PP^2}(-1)^6
 \longrightarrow i_*\beta\longrightarrow0.
\]
Its determinant is a nonzero scalar multiple of \(f_C\). Put
\[
 q_i(x)=x^tM_ix,\qquad S_\beta=V(q_0,q_1,q_2)\subset\PP^5.
\]
Every matrix in the net has rank at least five: a matrix of rank
at most four has zero adjugate and would give a singular point of
the determinant curve.

If the three differentials were dependent at \(p\in S_\beta\),
then \(M(\lambda)p=0\) for some \([\lambda]\). This matrix has
corank one, and its adjugate is a nonzero scalar multiple of
\(pp^t\). The derivatives of \(\det M(\lambda)\) are then scalar
multiples of \(p^tM_ip=q_i(p)=0\). This again contradicts the
smoothness of \(C\). The base locus is therefore a smooth complete
intersection of three quadrics, hence a K3 surface of degree eight.
\end{proof}

We recall the quadric-reduction construction underlying \cite[Example 4.10(a) and Corollary 4.11]{BKK}, giving the details needed under our hypotheses.

\begin{lemma}\label{lem:rational-reduction}
Under the hypotheses of Lemma~\ref{lem:net}, there is
a flat quadric surface bundle \(r:Z\to\PP^2\)
with integral total space \(Z,\) discriminant
double cover \(B\) and Clifford class \(\mathfrak b_B(\beta)\).
Moreover, \(Z\) is rational.
\end{lemma}
\begin{proof}
Choose \(p\in S_\beta\), and take coordinates with
\(p=[1:0:\cdots:0]\). Write
\[
 q_i(x_0,z)=x_0\ell_i(z)+r_i(z),\qquad z\in\CC^5,
\]
where \(\ell_i\) is linear and \(r_i\) is quadratic. Smoothness at
\(p\) implies that \(\ell_0,\ell_1,\ell_2\) are independent.
Thus \(\sum\lambda_i\ell_i\ne0\) for every \([\lambda]\), so
\(p\) is a smooth point of every quadric in the net. The constant
isotropic section is therefore regular.
Reduction of the net at \(p\) is the bundle
\begin{equation}\label{eq:reduced-net}
 Z=\left\{(\lambda,z)\in\PP^2\times\PP^4:
       \sum_{i=0}^2\lambda_i\ell_i(z)=
       \sum_{i=0}^2\lambda_i r_i(z)=0\right\}.
\end{equation}
The first equation gives a \(\PP^3\)-bundle over \(\PP^2\).
The second is a relative quadric. Reduction removes a hyperbolic
plane from a form of rank at least five, so every fiber has rank
at least three. Thus the fibers are geometrically integral quadric
surfaces, and the relative equation gives flatness. The total space
is integral. The Clifford class of the symmetric resolution is
\(\mathfrak b_B(\beta)\), and quadric reduction preserves it
\cite[Theorem~2.3 and Proposition~3.2]{BKK}.

Set
\[
 \ell(z)=(\ell_0(z),\ell_1(z),\ell_2(z)),\qquad
 r(z)=(r_0(z),r_1(z),r_2(z)).
\]
When these vectors are independent, the two equations in
\eqref{eq:reduced-net} determine a unique parameter
\[
 [\lambda]=[\ell(z)\times r(z)].
\]
This gives a rational inverse to the projection \(Z\to\PP^4\).
The open set is nonempty. Otherwise all
\(\ell_i r_j-\ell_j r_i\) would vanish identically. Since the
independent linear forms are pairwise coprime, this would imply
\(r_i=\ell_i L\) for a common linear form \(L\). Then
\(q_i=(x_0+L)\ell_i\), contradicting Lemma~\ref{lem:net}.
Hence \(Z\) is birational to \(\PP^4\).
\end{proof}

\begin{corollary}\label{cor:theta-rational}
Let \(q:W\to\PP^2\) be a quadric surface bundle with smooth
discriminant sextic \(C\) and double cover \(B\). If 
\(\alpha_W=\mathfrak b_B(\beta)\) for an ineffective even theta 
characteristic \(\beta\) on \(C\), where \(\mathfrak b_B\) is the 
homomorphism in (\ref{Sumit}), then \(W\) is rational.
\end{corollary}
\begin{proof}
Lemma~\ref{lem:rational-reduction} gives a quadric surface bundle 
\(Z\to\mathbf P^2,\) whose total space is rational, and whose
discriminant double cover and Clifford class agree with those of \(W.\) By
Proposition~\ref{quad}, \(W\) and \(Z\) are
birational over \(\PP^2\).
\end{proof}

Combining the vanishing Brauer-class calculation of 
Galluzzi–van Geemen \cite[Proposition 5.1.4 and its proof]{GvGbrauer} 
with the rational 
quadric-bundle construction of Bini–Kapustka–Kapustka \cite[Example 4.10(a) 
and Corollary 4.11]{BKK}, we obtain the following proposition.

\begin{proposition}\label{prop:theta-class}
Let \(X\) have full algebraic lattice \(M_{\tau,n}\), with
\(\tau=0,4\) and \(n\) odd. Suppose that its plane projection has
smooth discriminant sextic \(C\), with double cover \(B\).
There is an even theta characteristic \(\beta\) on \(C\) such that
\[
 \mathfrak b_B(\beta)=\alpha_X.
\]
It is ineffective, and \(X\) is rational.
\end{proposition}

\begin{proof}
By \cite[Theorem~4.3.1, Sections~4.10 and~4.14]{GvGbrauer},
there is an integral basis 
$(h,\kappa)$ of $\operatorname{NS}(B)$ in which
\begin{equation}\label{sad}
\operatorname{Gram}(h,\kappa)=
\begin{cases}
\begin{pmatrix}
2&0\\
0&-2n
\end{pmatrix},
& \tau=0,\\[6pt]
\begin{pmatrix}
2&0\\
0&6-2n
\end{pmatrix},
& \tau=4.
\end{cases}
\end{equation}
In particular, $h\cdot\operatorname{NS}(B)\subset 2\mathbf Z$.

We recall the calculation of the specialized Brauer class in
\cite[proof of Proposition~5.1.4]{GvGbrauer}. Deform the pair
\((X,P)\) in the parameter space of smooth cubics containing a
fixed plane. The additional algebraic classes define a countable
union of proper Hodge loci. A smooth local one-parameter family
through \((X,P)\) can therefore be chosen with very general fiber
of algebraic rank two. After shrinking the base, its discriminant
double covers are smooth, and the very general double plane
\(B_t\) has \(\NS(B_t)=\ZZ h\).

Choose a very general parameter $t$ and trivialize the
integral cohomology local system. Using \eqref{sad}, we may write
\[
\operatorname{NS}(B_t)=\mathbf Z h,
\qquad
\operatorname{NS}(B)=\mathbf Z h\oplus\mathbf Z\kappa
\]
inside a fixed K3 lattice $\Lambda$. Consequently,
\begin{equation}\label{happy}
T(B)
=
(\mathbf Z h+\mathbf Z\kappa)^\perp_\Lambda
=
T(B_t)\cap\kappa^\perp.
\end{equation}

The vanishing Brauer class \(\delta\) is represented by
\(\kappa/2\); see \cite[Proposition~2.1.2]{GvGbrauer}.
Write \(\alpha_t=\alpha_{X_t}\) for the Clifford class of the nearby cubic.
If \(\theta_t\) is the original odd theta characteristic,
the cited calculation supplies \(\eta_t\in\Pic(C_t)[2]\) with
\[
 \beta_t=\theta_t\otimes\eta_t,\qquad
 \mathfrak b_{B_t}(\beta_t)=\alpha_t+\delta.
\]

For these parameters, the theta characteristic \(\beta_t\) is even because 
\(n\) is odd.
After a finite base change if needed, it extends in the smooth
family to an even theta characteristic \(\beta\) on \(C\).

Under the description of the K3 Brauer group as homomorphisms
on its transcendental lattice, the restriction calculation in
\cite[proof of Proposition~5.1.4]{GvGbrauer} gives, for
\(x\in T(B)\),
\[
 \mathfrak b_B(\beta)(x)
   =\alpha_t(x)+\frac{\kappa\cdot x}{2}\pmod\ZZ
   =\alpha_X(x).
\]
The correction term \(\frac{\kappa \cdot x}{2}\) vanishes. Equality on \(T(B)\)
is equality in \(\Br(B)[2]\). Lemma~\ref{lem:ineffective} now
gives \(h^0(C,\beta)=0\), and
Corollary~\ref{cor:theta-rational} proves rationality.
\end{proof}

\subsection{Rationality through Verra fourfolds}

The use of the two quadric surface fibrations of a Verra fourfold to 
establish rationality has precedents in 
\cite[Proposition 6.3 and Corollary 6.4]{CKKM}. 
Here the existence of a section is deduced from the discriminant bilinear 
pairing of the algebraic lattice.

\begin{proposition}\label{prop:second-projection}
Let \(X\) have full algebraic lattice \(M_{2,n}\), with \(n\) odd.
Suppose there exists a smooth Verra fourfold \(V\) birational to \(X\), 
with one projection birational over \(\mathbf P^2\) to \(q_X\), 
and that the common Clifford data induce an integral Hodge isometry
\(T(V)\simeq T(X)\). Then the second Verra projection has a
rational section. Consequently, \(V\) and \(X\) are rational.
\end{proposition}
\begin{proof}
Label the projections \(p_1,p_2\) so that \(p_1\) is birational
over the base to the plane projection of \(X\), and put
\(N=A(V)\). This isometry, together with unimodularity and 
the ranks and signatures of the 
ambient cohomology lattices, gives
\[
 \rk N=4,\qquad \operatorname{sign}(N)=(3,1),\qquad
 \det N=-(16n-12).
\]
By Lemma~\ref{lem:discriminant},
\begin{equation}\label{eq:pairing-equality}
 (D(N),b_N)\simeq(D(M_{2,n}),b_{M_{2,n}}).
\end{equation}

The primitive lattice \eqref{eq:Verra-lattice} can be extended to
an integral basis \(f_1,f_2,e,w\) of \(N\). Write
\[
 \Gram(f_1,f_2,e,w)=
 \begin{pmatrix}
 0&2&0&a\\2&0&0&b\\0&0&2&c\\a&b&c&d
 \end{pmatrix}.
\]
Its determinant is \(4(2ab+c^2-2d)\).

Suppose that neither projection has a rational section.
Proposition~\ref{death} implies that \(a\) and \(b\)
are even. Replace \(w\) by
\[
 w-\frac a2 f_2-\frac b2 f_1.
\]
Then \(a=b=0\), and the determinant identity becomes
\[
 2d-c^2=4n-3.
\]
In particular \(c\) is odd. Replacing \(w\) by
\(w-(c-1)e/2\) makes \(e\cdot w=1\); the determinant then gives
\(w^2=2n-1\). Hence
\begin{equation}\label{eq:Verra-no-sections}
 N\simeq U(2)\oplus
 \begin{pmatrix}2&1\\1&2n-1\end{pmatrix}.
\end{equation}
The second summand has odd determinant \(4n-3\). Every element
of \(D(N)[2]\) therefore comes from \(U(2)\), with representatives
\[
 0,\qquad\frac{f_1}{2},\qquad\frac{f_2}{2},\qquad
 \frac{f_1+f_2}{2}.
\]
Their squares are \(0,0,0,1\), respectively. Thus
\[
 b_N(\xi,\xi)=0\quad\text{in }\QQ/\ZZ
 \quad\text{for every }\xi\in D(N)[2].
\]

On the cubic side, \(R/2\in M_{2,n}^\vee\), since its pairings
with \(H,P,R\) are \(0,1,n\). Its discriminant class has order
two and self-pairing
\[
 b_{M_{2,n}}([R/2],[R/2])=\frac n2=\frac12\pmod\ZZ,
\]
because \(n\) is odd. This contradicts
\eqref{eq:pairing-equality}. At least one projection has a rational section.

By \eqref{eq:Q-pairings}, the residual quadric pairs evenly with
every element of \(A(X)=M_{2,n}\). Thus the plane projection
has no rational section by \eqref{life}. Its generic
quadric is birational to that of \(p_1\), so \(p_1\) has no
rational section either. Therefore \(p_2\) has a rational section. 
Let \(K=\mathbf C(\mathbf P^2)\) denote the function field of the base of \(p_2\).
Its generic quadric is \(K\)-rational, and
\[\CC(X)\simeq\CC(V)\simeq K(u,v)\simeq\CC(x,y,u,v).\]
This proves the proposition.
\end{proof}

\begin{corollary}\label{cor:general-rational}
A very general member of \(\cC_{\tau,n}\) is rational whenever
\(M_{\tau,n}\) is admissible.
\end{corollary}
\begin{proof}
For \(\tau=1,3\), use \eqref{eq:Q-pairings} and Hassett's section
criterion. For \(\tau=0,4\) and odd \(n\), the general plane
projection has smooth discriminant, and
Proposition~\ref{prop:theta-class} applies. The exceptional
multiple-plane cases in the rank-three classification do not
occur in these two series.

For \(\tau=2\) and odd \(n\),
\cite[Theorem~0.1]{GvGVerra} supplies the general smooth compatible
Verra model. By \cite[Proposition 1.3(b) and Section 1.7]{GvGVerra}, the common discriminant K3 surface and Clifford invariant give an integral Hodge isometry \(T(V)\simeq T(X)\), as explained in Section~2.3. Proposition~3.7 therefore applies.
\end{proof}

\subsection{Smooth specialization and higher algebraic rank}
The following lemma is an immediate consequence of
\cite[Theorem~1 and the following paragraph in the introduction]{KT}.

\begin{lemma}
Let \(\pi:\mathcal X\to B\) be a smooth projective family
over an irreducible complex algebraic variety.
If its very general fiber is rational, then every fiber
over a closed point of \(B\) is rational.
\end{lemma}

Combining Corollary 3.8 with the specialization theorem of Kontsevich–Tschinkel 
and the irreducibility established by Yang–Yu, we obtain the following theorem.

\begin{theorem}\label{good}
If \(M_{\tau,n}\) is admissible, every smooth cubic fourfold
in \(\mathcal C_{\tau,n}\) is rational.
\end{theorem}

\begin{proof}
Let
\[
U\subset
\mathbf P\bigl(H^0(\mathbf P^5,\mathcal O_{\mathbf P^5}(3))\bigr)
\]
be the open subset parametrizing smooth cubic fourfolds, and let
\[
q:U\longrightarrow\mathcal C
\]
be the geometric GIT quotient by
\(G=\operatorname{PGL}(6,\mathbf C)\); see \cite[Section~3]{YY}.
By \cite[Corollary~8.14]{YY}, the locus
\(\mathcal C_{\tau,n}\) is closed and irreducible. Put
\[
B=\bigl(q^{-1}(\mathcal C_{\tau,n})\bigr)_{\mathrm{red}}.
\]

The variety \(B\) is irreducible. Indeed, since \(G\) is
connected, each irreducible component of \(B\) is
\(G\)-invariant. Its image under the GIT quotient is closed.
These finitely many images cover the irreducible variety
\(\mathcal C_{\tau,n}\), so one component maps onto
\(\mathcal C_{\tau,n}\). Each fiber of \(q\) is a single
\(G\)-orbit. The component mapped onto
\(\mathcal C_{\tau,n}\) therefore contains every fiber over
this locus and hence equals \(B\).

Restricting the universal family of smooth cubic fourfolds
over \(U\) gives a smooth projective family
\[
\pi:\mathcal X_B\longrightarrow B.
\]
By Corollary~3.8, a very general member of
\(\mathcal C_{\tau,n}\) is rational. Since
\(q|_B:B\to\mathcal C_{\tau,n}\) is surjective, the inverse
image of each proper closed exceptional subset is again
proper and closed. Thus the very general fiber of
\(\pi\) is rational.

Lemma~3.9 now implies that \(\mathcal X_b\) is rational
for every closed point \(b\in B\). Every cubic fourfold
in \(\mathcal C_{\tau,n}\) is represented by such a point,
which proves the theorem. In particular, the conclusion
also holds when the algebraic rank increases or the
discriminant of the plane projection becomes singular:
the family of cubic fourfolds itself remains smooth.
\end{proof}

The following lemma is a consequence of the
construction in \cite[proof of Proposition~7.3]{YY}; see also
\cite[Corollary~7.5]{YY}. We apply this construction to the plane
marking \(K_8=\langle H,P\rangle\) and an admissible rank-two
marking \(K_d\subset A(X)\). Since \(d\ne8\), the saturation 
of \(K_8+K_d\) inside \(A(X)\) is rank three. It contains \(K_d\)
primitively and is therefore admissible in the sense of
\cite[Definition~8.1]{YY}. Its normal form \(M_{\tau,n}\) is given
by \cite[Corollary~8.14]{YY}.

\begin{lemma}\label{bad}
Let \(X\in\mathcal C_8\), and suppose that \(A(X)\) is
admissible. Then \(A(X)\) contains a primitive admissible
rank-three lattice containing \(H\) and the class of a plane.
\end{lemma}

\begin{proposition}\label{prop:sufficiency}
Let \(X\in\cC_8\). If \(A(X)\) is admissible, then \(X\) is rational.
\end{proposition}
\begin{proof}
By Lemma~\ref{bad}, \(X\) belongs to a locus
with a primitive admissible \(M_{\tau,n}\)-marking.
Theorem~\ref{good} applies.
\end{proof}

\subsection{Necessity and the proof of the main theorem}
We use the transcendental-lattice conclusion stated in 
\cite[§3.3, Proposition 26 and the proof of Theorem 27]{GuereIntegral}. 
Our argument requires an integral Hodge isometry to \(T(X)\simeq T(S)(-1)\); 
it does not require an integral isometry between 
the full primitive cohomology of \(X\) 
and the full second cohomology of a K3 surface.

\begin{proposition}[J. Guéré]\label{god}
Let \(X\) be a rational smooth cubic fourfold.
Then there is a projective K3 surface \(S\) and an integral
Hodge isometry \[T(X)\simeq T(S)(-1).\]
In particular, \(A(X)\) is admissible.
\end{proposition}

\begin{proof}
Choose a weak factorization of a birational map
\(X \dashrightarrow \mathbf P^4\),
\[
X=X_0 \dashleftarrow\dashrightarrow X_1
\dashleftarrow\dashrightarrow \cdots
\dashleftarrow\dashrightarrow X_m=\mathbf P^4.
\]
Work over an auxiliary number field \(K\), enlarging it as required
for the comparison of evaluation maps along this factorization.

Choose a general evaluation \(\mathrm{ev}\), in the open set where
the spectrum has maximal cardinality. Let
\(E^X_{\mathrm{ev},\alpha}\) be the atom containing \(T(X)\)
after extension of coefficients. The last paragraph of
\cite[proof of Proposition~62]{GuereCubic} gives
\[
\nu^X_{\mathrm{ev},\alpha}=1,\qquad
\nu'{}^X_{\mathrm{ev},\alpha}=0,\qquad
\gamma^X_{\mathrm{ev},\alpha}\leq 1,
\]
where we use the notation of
\cite[Definition~28]{GuereCubic}.

Follow this atom along the factorization, using
\cite[Remark~32, Corollary~43, and Remark~48]{GuereCubic}.
At a blow-up, choose general evaluations on the added center
summands with disjoint spectra. At a blow-down, the transported
atom belongs to a unique summand. These comparisons preserve the
displayed numerical invariants. Since \(\mathbf P^4\), points,
and curves have no contribution to \(\nu\), the atom must pass
into a surface center \(Y\).

The integral decomposition in
\cite[Sections~3.2--3.3 and Propositions~25--26]{GuereIntegral}
is compatible with these atom comparisons. This compatibility
follows from
\cite[Lemmas~16--17 and Proposition~18]{GuereIntegral}.
Consequently, the integral transcendental summand of the chosen
atom identifies with \(T(Y)(-1)\) for this same surface center.
Since the chosen atom contains all of \(T(X)\), we obtain an
integral Hodge isometry
\[
T(X)\simeq T(Y)(-1).
\]

The induced evaluation on \(Y\) violates the strong property
\(\heartsuit\). By
\cite[Lemma~38]{GuereCubic}, the intrinsic strong property
fails as well. Taking \(K\) to contain \(\mathbf Q(i)\),
\cite[Proposition~60]{GuereCubic} gives failure of the
generic property. Hence
\cite[Proposition~59]{GuereCubic} shows that the smooth
projective minimal model \(S\) of \(Y\) satisfies
\[
c_1(K_S)=0 \quad\text{in }H^2(S,\mathbf Q),
\qquad
h^1(S,\mathcal O_S)=0.
\]

Point blow-ups preserve the integral transcendental lattice,
as recalled in Section~2.4. Thus
\[
T(X)\simeq T(Y)(-1)\simeq T(S)(-1).
\]
In particular,
\[
h^2(S,\mathcal O_S)=h^{2,0}(Y)=h^{3,1}(X)=1.
\]
A nonzero section of \(K_S\) has an effective zero divisor
numerically equivalent to zero. Such a divisor must vanish,
since every nonzero effective divisor has positive intersection
with an ample class. Hence \(K_S\simeq\mathcal O_S\), and
\(h^1(S,\mathcal O_S)=0\) shows that \(S\) is a K3 surface.

The transcendental-lattice criterion recalled in Section~2.1
now implies that \(A(X)\) is admissible.
\end{proof}

\begin{proof}[Proof of Theorem~1.1]
If \(X\) is rational, Proposition~3.13 shows that
\(A(X)\) is admissible. Conversely, if \(A(X)\) is admissible,
Proposition~3.12 shows that \(X\) is rational.
\end{proof}

\begin{remark}
By definition, \(A(X)\) is admissible precisely when it
contains a primitive rank-two marking \(K_d\) whose
discriminant \(d\) satisfies (2). This is equivalent to
\(X\in\mathcal C_d\) for some such \(d\). Restricting to
\(\mathcal C_8\) gives the asserted description
\[
\bigcup_{\substack{d>6,\ d\equiv0,2\;(\mathrm{mod}\,6)\\
d\text{ satisfies }(2)}}
(\mathcal C_8\cap\mathcal C_d).
\]
of the rational locus.
\end{remark}

\begin{proof}[Proof of Corollary~1.2]
Assume that \(A(X)=M_{\tau,n}\), with the parameters
normalized as in the statement. By Proposition~3.1,
this lattice is admissible if and only if
\[
\tau\in\{1,3\}
\quad\text{or}\quad
\bigl(\tau\in\{0,2,4\}\text{ and }n\text{ is odd}\bigr).
\]
The equivalent with rationality follows from
Theorem~1.1. For each admissible pair \((\tau,n)\),
Theorem~3.10 proves rationality for every smooth cubic
with a primitive \(M_{\tau,n}\)-marking.
\end{proof}

\begin{remark}
The nonrationality assertion in Corollary~1.2 assumes
\(A(X)=M_{\tau,n}\). If \(X\) merely admits a primitive
\(M_{\tau,n}\)-marking, inadmissibility of \(M_{\tau,n}\)
does not imply inadmissibility of \(A(X)\): the larger
lattice may contain an admissible rank-two marking.
In particular, a smooth specialization may acquire such
a marking and hence become rational. The nonrationality
conclusion must therefore be tested on the full algebraic
lattice. An admissible rank-three marking remains
sufficient for rationality when the algebraic rank
increases, by Theorem~3.10.
\end{remark}

\end{document}